\documentclass{amsart}
\usepackage{amsmath}

\theoremstyle{plain}

\newtheorem*{ther}{Theorem}

\theoremstyle{definition}

\def \CPb {\overline{\mathbf{CP}}^{\,2}}

\def \Z {\mathbf{Z}}
\def \Sig{\Sigma}

\def \la {\langle}
\def \ra {\rangle}
\def \a {\alpha}

\def \s {\sigma}
\def \t {\tau}

\def \bd {\partial}
\def \x {\times}
\def \- {\setminus}
\def \C {\subset}

\def\spinc{spin$^{\text{c}}$} 
\def\ct{\mathcal{T}}

\def \ssw {\text{SW}}
\def \sw {\mathcal{SW}}

\def \DD {\Delta}

\begin{document}

\baselineskip.475cm

\title{ Surfaces in 4-Manifolds: Addendum}
\author[Ronald Fintushel]{Ronald Fintushel}
\address{Department of Mathematics, Michigan State University \newline
\hspace*{.375in}East Lansing, Michigan 48824}
\email{\rm{ronfint@math.msu.edu}}
\thanks{The first author was partially supported NSF Grant DMS0305818
and the second author by NSF Grant DMS0204041}
\author[Ronald J. Stern]{Ronald J. Stern}
\address{Department of Mathematics, University of California \newline
\hspace*{.375in}Irvine,  California 92697}
\email{\rm{rstern@math.uci.edu}}

\begin{abstract} In this note we fill a gap in the proof of the main theorem (Theorem 1.2)\/ of our paper
{\em Surfaces in 4-manifolds}, Math. Res. Letters {$\mathbf{4}$} (1997), 907--914.
\smallskip
\end{abstract}

\maketitle

Let $\Sig$ be a smoothly embedded surface in a simply connected smooth 4-manifold, and assume that $\Sig$ has nonnegative self-intersection number $n$ and  satisfies 
$\pi_1(X\setminus \Sig)=0$. Given a knot $K$ in $S^3$ and a nontrivial loop $C$ on 
$\Sig$, one can perform `rim surgery' by choosing a trivialization $C\x D^2$ of the normal bundle of $\Sig$ restricted over $C$ and performing knot surgery \cite{KL4M} on the rim torus $C\x \bd D^2$ in $X$. This operation gives a new surface 
$\Sig_{K,C}\C X$. 

Let $\Sig_n$ be the surface of self-intersection $0$ in $X_n=X\# n\CPb$ obtained by blowing up at $n$ points of $\Sig$. 
In \cite{surfaces} we defined a collection of complex surfaces $Y_g$ containing standardly embedded surfaces $S_g$ of self-intersection $0$, and we called $(X,\Sig)$ a {\it SW-pair}\/ if the Seiberg-Witten invariant $\sw_{X_n\#_{\Sig_n=S_g} Y_g}\ne 0$. The main theorem of \cite{surfaces} states that for any $\ssw$-pair $(X,\Sig)$, if
$K_1$ and $K_2$ are two knots in $S^3$ and if there is a diffeomorphism of pairs
$(X, \Sig_{K_1})\to (X,\Sig_{K_2})$, then  $\DD_{K_1}(t)=\DD_{K_2}(t)$.
To prove this theorem, one identifies $X_n\#_{\Sig_{n,K}=S_g} Y_g$ with the result of knot surgery on the fiber sum $(X_n\#_{\Sig_{n}=S_g} Y_g)_K$, and then uses
\cite{KL4M} to calculate that $\sw_{X_n\#_{\Sig_{n,K}=S_g} Y_g}=\sw_{X_n\#_{\Sig_n=S_g} Y_g}\cdot\DD_{K}(t^2)$. 

However, some issues were not properly addressed. In particular, the fiber sum construction $X_n\#_{\Sig_n=S_g} Y_g$ needs both a fixed identification of the surface $\Sig_n$ with $S_g$ as well as a choice of framing. The statement that $X_n\#_{\Sig_{n,K}=S_g} Y_g$ is diffeomorphic to $(X_n\#_{\Sig_{n}=S_g} Y_g)_K$ assumes that the identification and framing for the first fiber sum is induced from the choices made for $X_n\#_{\Sig_{n}=S_g} Y_g$. Since $(Y_1,S_1) =(E(1),{\text{fiber}})$, and since the complement of a fiber in $E(1)$ has big diffeomorphism group, there is no problem when $g=1$. 

Although it is possible to repair this problem, the construction of monopole Floer homology by Kronheimer and Mrowka \cite{book} (cf.\,{\cite{KMOS}})  gives us a more satisfactory method of dealing with this situation. The point of fiber-summing $X_n$ to $Y_g$ was to exhibit the effect of rim surgery on the relative Seiberg-Witten invariant of $X_n\- \Sig$. The \spinc-structures on $\Sig\x S^1$ are in $1-1$ correspondence with elements of $H^2(\Sig\x S^1;\Z)\cong H^2(\Sig;\Z) \oplus H^1(\Sig;\Z) \cong \Z\oplus H^1(\Sig;\Z)$; but any
\spinc-structure not pulled back from $\Sig$ has a trivial Floer homology group. Thus one only needs to consider the \spinc-structures $s_k$ corresponding to $(k,0)$. This \spinc-structure satisfies $\la c_1(s_k),[\Sig]\ra =2k$.
The Floer homology group $HM(\Sig\x S^1; s_k)$ is trivial for $|k|\ge g$, the genus of $\Sig$. We are interested in $HM(\Sig\x S^1; s_{g-1})\cong \Z$. 

For convenience we now assume that $\Sig\cdot\Sig =0$; so $(X_n,\Sig_n)= (X,\Sig)$. The relative Seiberg-Witten invariant $\sw_{X, \Sig}$ assigns to each \spinc-structure $\t$ on $X\- N(\Sig)$ an element in $HM(\Sig\x S^1; \s)$ where $\s$ is the \spinc-structure on 
$\Sig\x S^1 = \bd N(\Sig)$ obtained by restricting $\s$. Let $\ct$ be the collection of \spinc-structures $\t$ on $X\- N(\Sig)$ whose restriction to $\bd N(\Sig)$ is $\pm s_{g-1}$. This gives rise to a well-defined Seiberg-Witten invariant $\ssw^{\ct}_{X,\Sig}:\ct\to\Z$.  In the usual way,one obtains a Laurent polynomial $\sw^{\ct}_{X,\Sig}$ with variables from $A=\{\a\in H^2(X\- \Sig;\Z)\mid \a|_{\Sig\x S^1}=s_{g-1}\}$. This is an invariant in the sense that a diffeomorphism $f: (X,\Sig)\to (X',\Sig')$ induces $f^*: A'\to A$ and sends $\sw^{\ct}_{X',\Sig'}$ to $\sw^{\ct}_{X,\Sig}$. 
Note that $H^2(X\- \Sig;\Z)$ has a summand isomorphic to $H^3(X,X\- \Sig;\Z)\cong H_1(\Sig;\Z)\cong R$, the subgroup of $H_2(\Sig\x S^1;\Z)$ generated by the rim tori of $\Sig$. There is a canonical identification of $H^2(X\- \Sig;\Z)$ and $H^2(X\- \Sig_K;\Z)$ which identifies $A=A_{(X,\Sig)}$ with $A_{(X,\Sig_K)}$.

\begin{ther} Let $\Sig$ be a smoothly embedded genus $g>0$ surface in a simply connected smooth 4-manifold $X$, and assume that $\Sig$ has self-intersection number $0$ and  satisfies $\pi_1(X\setminus \Sig)=0$. Assume that the relative Seiberg-Witten invariant $\sw^{\ct}_{X,\Sig}\ne 0$. If
$K_1$ and $K_2$ are two knots in $S^3$ and if there is a diffeomorphism of pairs
$f: (X,\Sig_{K_1})\to (X,\Sig_{K_2})$, then  the set of coefficients (with multiplicities) of $\DD_{K_1}(t)$ must be equal to that of $\DD_{K_2}(t)$.
\end{ther}
\begin{proof} The proof of the knot surgery theorem \cite{KL4M} applies in this situation to show that $\sw^{\ct}_{X,\Sig_{K_i}} = \sw^{\ct}_{X,\Sig}\cdot \DD_{K_i}(r_i^2)$ where $r_i$ is the element of $R$ belonging to the rim torus on which rim surgery was done.
  The theorem follows because the coefficients of $\sw^{\ct}_{X,\Sig}\cdot \DD_{K_1}({r_1}^2)$ must be precisely equal to those of $\sw^{\ct}_{X,\Sig}\cdot \DD_{K_2}({r_2}^2)$.
\end{proof}

\noindent {\it{Remarks.}} \ 1. That the hypothesis $\sw^{\ct}_{X,\Sig}\ne 0$ is weaker than the hypothesis 
$\sw_{X\#_{\Sig=S_g}Y_g}\ne 0$ of \cite{surfaces} follows from the gluing formula \cite{book} because $Y_g$ is a complex surface, and its canonical class $K$ gives a basic class for which $K\cdot S_g = 2g-2$.  

\noindent 2. The conclusion of the theorem is slightly weaker than that claimed in \cite{surfaces}. The authors currently see no way around this, although it is conceivable that the hypothesis of the theorem implies that $\DD_{K_1}(t)=\DD_{K_2}(t)$. For example, as Stefano Vidussi has pointed out, this is true in case $\sw^{\ct}_{X,\Sig}$ is an irreducible (over $\Z$) Laurent polynomial whose support in $A$ has rank at least 2.

\noindent 3. The authors wish to thank Stefano Vidussi and Danny Ruberman for helpful discussions.

\end{document}